# PLURISUBHARMONIC FUNCTIONS AND SUBELLIPTICITY OF THE $\bar{\partial}$-NEUMANN PROBLEM ON NONSMOOTH DOMAINS

Emil J. Straube*

*1. Introduction*

Recently, Michel and Shaw ([15]) and Henkin, Iordan, and Kohn [11] proved 1/2-subellipticity of the $\bar{\partial}$-Neumann problem on domains with piecewise smooth strictly pseudoconvex boundaries. (A more precise description of their results is in section 5.) In this paper, we discuss the situation when the boundary is piecewise smooth of finite type. It is shown that near such a boundary point, a subelliptic estimate holds for some $\varepsilon > 0$ (Theorem 2). We use Catlin's result on the existence of families of plurisubharmonic functions whose Hessians satisfy good lower bounds near finite type points ([5]). What is new is that just as in the smooth case ([5], Theorem 2.2), the existence of (families of) plurisubharmonic functions with algebraic growth of the Hessian gives a subelliptic estimate when the boundary is only Lipschitz (Theorem 1, Remark 1). Since the boundary is only Lipschitz, the argument is necessarily different from [5]. By using the fact that the (formal) complex Laplacian acts diagonally on forms as a multiple of the real Laplacian, we first reduce the problem to the case of forms with harmonic coefficients. Following [4], (2.3), [5], Theorem 2.1 for this part of the argument, we obtain from the existence of the good plurisubharmonic functions that $\|\bar{\partial}u\|_0 + \|\bar{\partial}^*u\|_0$ dominates the $\mathcal{L}^2$-norm of the form $u$ weighted by an inverse power of the boundary distance. For harmonic functions, however, it is well understood that this weighted norm controls a Sobolev norm. Due to the minimal smoothness assumptions, an additional complication arises. To obtain the domination of the weighted $\mathcal{L}^2$-norm of $u$, one needs to work first on smooth subdomains. The problem that the restrictions of the forms to the subdomains are not, in general, in the domain of $\bar{\partial}^*$ there, is handled by an argument involving the $\bar{\partial}$-Neumann operator of the subdomains to modify these restrictions. This regularization procedure is useful in related contexts. For example, it allows to extend Catlin's classical result on compactness of the $\bar{\partial}$-Neumann operator in the presence of plurisubharmonic functions with large Hessians

---

* Supported in part by NSF grant DMS-9500916.

([4], Theorem 1) to the case where no boundary smoothness at all is assumed (Corollary 3). This simplifies and improves existing compactness results in the nonsmooth case.

*Acknowledgement:* The author thanks Mei-Chi Shaw for catching an error in an earlier version and for a stimulating conversation on the subject of this paper.

## 2. Statement of results

In this paper, we will be mainly concerned with bounded domains in $C^n$ with Lipschitz boundary. This means that the boundary is locally the graph of a Lipschitz function. By $W^s(\Omega)$ for $s \geq 0$ we denote the $\mathcal{L}^2$-Sobolev space of order $s$ on $\Omega$, i.e. the restrictions of functions in $W^s(C^n)$ to $\Omega$, with the quotient norm. For $s \in N$, this gives the space of functions which together with their derivatives of order up to $s$ are square integrable. (For a discussion of Sobolev spaces on Lipschitz domains, see for example [13], section 2.) It will be convenient to adopt the same definition of $W^s(\Omega)$ also when $\Omega$ does not have Lipschitz boundary. $\delta(z) = \delta_\Omega(z)$ denotes the distance to the boundary. We let $d^c = i(\bar{\partial} - \partial)$, so that $dd^c = 2i\partial\bar{\partial}$.

For background on the $\bar{\partial}$-Neumann problem, we refer the reader to [9], [6]. The basic $\mathcal{L}^2$-theory is conveniently summarized in [3], Proposition 3, and [16], section 2. For a survey of the $\mathcal{L}^2$-Sobolev theory, see [1]. For $1 \leq q \leq n$, the $\bar{\partial}$-Neumann operator $N_q$ is the inverse of the complex Laplacian $\bar{\partial}^*\bar{\partial} + \bar{\partial}\bar{\partial}^*$ on $(0, q)$-forms.

Let $P$ be a boundary point of the bounded pseudoconvex domain $\Omega$. A subelliptic estimate of order $\varepsilon$ is said to hold near $P$ if there are a neighborhood $U$ of $P$, $\varepsilon > 0$, and a constant $C > 0$, such that

$$\|u\|_{\varepsilon, U \cap \Omega} \leq C(\|\bar{\partial}u\|_0 + \|\bar{\partial}^*u\|_0), \tag{1}$$

for all $u \in \text{dom}(\bar{\partial}) \cap \text{dom}(\bar{\partial}^*)$. For a recent survey of subellipticity in the case of smooth domains, see [7].

Our first theorem concerns subellipticity near a boundary point with good plurisubharmonic functions.

**Theorem 1.** *Let $\Omega$ be a bounded pseudoconvex domain, in $C^n$, $P$ a boundary point. Assume that the boundary is the graph of a Lipschitz function near $P$, and that there* exists a bounded plurisubharmonic function $\lambda$ in the intersection of a neighborhood $U$ of $P$ with $\Omega$, such that

$$dd^c\lambda \geq \text{const. } \delta^{-2\varepsilon} dd^c |z|^2 \text{ in } U \cap \Omega \qquad (2)$$

as currents, with $0 \leq \varepsilon < 1/2$. Then there is a neighborhood $V$ of $P$ and a constant $C > 0$, such that

(i)
$$\|u\|_{\varepsilon,\Omega \cap V} \leq C(\|\bar{\partial}u\|_0 + \|\bar{\partial}^*u\|_0),$$

$u \in \text{dom}(\bar{\partial}) \cap \text{dom}(\bar{\partial}^*) \subseteq \mathcal{L}^2_{(0,q)}(\Omega), 1 \leq q \leq n,$ \hfill (3)

and

(ii)
$$\|N_q u\|_{\varepsilon,\Omega \cap V} \leq C\|u\|_0,$$

$u \in \mathcal{L}^2_{(0,q)}(\Omega), 1 \leq q \leq n$.

For details about currents and their comparison, see for example [14], chapter 2. In view of the $\mathcal{L}^2$-continuity of $\bar{\partial} N_q$ and $\bar{\partial}^* N_q$, (ii) is of course a consequence of (i).

*Remark 1:* Compared with [5], Theorem 2.2, our assumption on the plurisubharmonic function with good growth of the Hessian is stronger in that we require *one* function whose Hessian blows up algebraically, rather than a (uniformly bounded) *family* of plurisubharmonic functions with suitable lower bounds on the Hessians. However, we will show in section 4 in the proof of Theorem 2 that given a family of plurisubharmonic functions as in [5], Theorem 2.2 for some $\varepsilon > 0$, one can easily construct a bounded plurisubharmonic function with the required growth of the Hessian as in Theorem 1 for every $\varepsilon' < \varepsilon$. From our point of view, the main point is the fact that the boundary is only assumed Lipschitz.

Theorem 1 applies to the situation where locally, near $P$, $\Omega$ is a transversal intersection of pseudoconvex domains whose boundaries contain $P$ and which are $C^\infty$-smooth and of finite type near $P$. To be precise, we assume that in suitable local coordinates $(z_1, \ldots, z_n)$, $z_j = x_j + iy_j$, $1 \leq j \leq n$, $\Omega$ is given as follows. There are $C^\infty$-functions $\rho_k(z) = y_n - h_k(z_1, \ldots, z_{n-1}, x_n)$, $1 \leq k \leq m$, so that, near $P$, $\Omega = \{\rho_k(z) < 0 \mid 1 \leq k \leq m\}$, the differentials $d\rho_k$, $1 \leq k \leq m$, are linearly independent over $R$, and the surfaces $\{\rho_k(z) = 0\}$ are pseudoconvex from the side $\{\rho_k < 0\}$ and of finite D'Angelo 1-type near $P$ ([6]). We say that $b\Omega$ is piecewise smooth of finite type near $P$.

**Theorem 2.** *Let $\Omega$ be a bounded pseudoconvex domain in $C^n, P$ a boundary point. Assume the boundary is piecewise smooth of finite type near $P$. Then the $\bar{\partial}$-Neumann problem is subelliptic near $P$, i.e. the conclusions of Theorem 1 hold.*

The regularizing procedure used in the proof of Theorem 1 has other related applications. For example, it immediately shows that Catlin's compactness theorem for the $\bar{\partial}$-Neumann problem ([4], Theorem 1) holds without any boundary smoothness assumption on the domain. The following is a simple corollary of part of the *proof* of Theorem 1.

**Corollary 3.** *Let $\Omega$ be a bounded pseudoconvex domain in $C^n$ with the property that for each positive number $M$ there exists a neighborhood $U$ of the boundary and a plurisubharmonic function $\lambda$ on $U \cap \Omega$, $0 \leq \lambda \leq 1$, with $dd^c\lambda \geq M dd^c|z|^2$ (as currents). Then the $\bar{\partial}$-Neumann operators on $\Omega$ are compact.*

In particular, the $\bar{\partial}$-Neumann operators are compact if the boundary of $\Omega$ (with no smoothness assumptions) is a $B$-regular set in the sense of Sibony ([17]). For the case of domains in $C^n$, Corollary 3 generalizes a recent result in [10].

*3. Proof of Theorem 1*

We adopt the customary convention on constants: $C$ may change its value as the argument progresses, but it will always be independent of the relevant forms. Let $u = \sum_J{}' u_J d\bar{z}_J \in \text{dom}(\bar{\partial}) \cap \text{dom}(\bar{\partial}^*)$. The prime in the summation indicates as usual that the summation is only over strictly increasing $q$-tuples $J$. We first reduce to the case where $u$ has harmonic coefficients. Denote by $\theta$ the formal adjoint of $\bar{\partial}$. Note that for $u \in \text{dom}(\bar{\partial}^*)$, $\bar{\partial}^* u = \theta u$. Since $\bar{\partial}\theta + \theta\bar{\partial}$ acts componentwise as a constant multiple of the (real) Laplacian, we have for all $J$

$$\|\Delta u_J\|_{-1} \leq C(\|\bar{\partial}u\|_0 + \|\bar{\partial}^*u\|_0). \tag{4}$$

Denote by $v_J$ the unique function in $W_0^1(\Omega)$ with $\Delta v_J = \Delta u_J$. (We use here that since $\Omega$ is bounded, $\Delta\colon W_0^1(\Omega) \to W^{-1}(\Omega)$ is an isomorphism, see for example [18], chapter 23.) Let $v := \sum_J{}' v_J d\bar{z}_J$. Then $v \in \text{dom}(\bar{\partial}) \cap \text{dom}(\bar{\partial}^*)$ (because all the $v_J$ are in $W_0^1(\Omega)$), and

$$\|v\|_1 \leq C\|\Delta u\|_{-1} \leq C(\|\bar{\partial}u\|_0 + \|\bar{\partial}^*u\|_0). \tag{5}$$

Then
$$u = (u - v) + v, \qquad (6)$$

$u - v \in \text{dom}(\bar{\partial}) \cap \text{dom}(\bar{\partial}^*)$ (since $u$ and $v$ are), and

$$\|\bar{\partial}(u-v)\|_0 + \|\bar{\partial}^*(u-v)\|_0 \leq C(\|\bar{\partial}u\|_0 + \|\bar{\partial}^*u\|_0 + \|v\|_1) \qquad (7)$$
$$\leq C(\|\bar{\partial}u\|_0 + \|\bar{\partial}^*u\|_0).$$

$u - v$ has harmonic coefficients, so once (3) is established for forms with harmonic coefficients, it will follow in full generality (in view of (5) and (7)).

We now assume that $u$ has harmonic coefficients. Let $V_0$ be a neighborhood of $P$ so that in suitable local coordinates $(z_1, \ldots, z_{n-1}, z_n) = (z', z_n)$, $\Omega \cap V_0$ is given by

$$\Omega \cap V_0 = \{z / |z'|^2 + x_n^2 < R^2, 0 < y_n < h(z', x_n)\} \qquad (8)$$

for some $R > 0$ and a (strictly positive) Lipschitz function $h$. We may assume that the $(z', x_n)$ coordinates of $P$ are $(0,0)$. Let $V_1$ be a neighborhood of $P$ that is relatively compact in $V_0$ and such that

$$\Omega \cap V_1 = \{z / |z'|^2 + x_n^2 < R^2/4, a < y_n < h(z', x_n)\} \qquad (9)$$

for some $a > 0$. Then (note that $\Omega \cap V_1$ has Lipschitz boundary)

$$\|u\|_{\varepsilon, V_1 \cap \Omega} \leq C(\|\delta_{V_1 \cap \Omega}^{1-\varepsilon} \nabla u\|_{0, V_1 \cap \Omega} + \|u\|_{0, V_1 \cap \Omega}), \qquad (10)$$

see [13], Theorem 4.2 (compare also Prop. 4.15, harmonicity is not needed at this step). Here, $\nabla u$ denotes the vector of all first order derivatives of (the components of) $u$. Note that (10) is a *genuine* estimate (as opposed to an *a priori* estimate): if the right-hand side is finite then $u$ *is* in $W_{(0,q)}^{\varepsilon}(V_1 \cap \Omega)$, and the inequality holds. The last term on the right-hand side of (10) is dominated by $\|u\|_0$, which is dominated by $\|\bar{\partial}u\|_0 + \|\bar{\partial}^*u\|_0$. Note that on $V_1 \cap \Omega$, $\delta_{V_1 \cap \Omega} \leq \delta_\Omega$, so that the first term on the right-hand side of (10) is dominated by $\|\delta^{1-\varepsilon} \nabla u\|_{0, V_1 \cap \Omega} = \|\delta^{-\varepsilon}(\delta \nabla u)\|_{0, V_1 \cap \Omega}$. Let $\varphi \in C_0^\infty(V_0)$ be identically 1 on a neighborhood of $\overline{V_1 \cap \Omega}$. We now exploit the harmonicity of $u$ via standard subaveraging properties of $|\nabla u|^2$ on balls centered at points $z$ of $V_1 \cap \Omega$ of radius comparable to $\delta(z)$, and contained in $\Omega$ and the set where $\varphi$ is one, to obtain that $\|\delta^{-\varepsilon}(\delta\nabla u)\|_{0,V_1\cap\Omega}$ is dominated by $\|\delta^{-\varepsilon}\varphi u\|_{0,V_0\cap\Omega}$ (see [8], proof of Lemma 1, [13], proof of Theorem 4.2). Summarizing, we have

$$\|u\|_{\varepsilon,V_1\cap\Omega} \leq C(\|\delta^{-\varepsilon}\varphi u\|_{0,V_0\cap\Omega} + \|\bar{\partial}u\|_0 + \|\bar{\partial}^*u\|_0). \tag{11}$$

We choose $V_0$ small enough so that $\overline{V_0\cap\Omega}$ is a compact subset of the neighborhood $U$ that appears in the hypothesis of the theorem. We may also assume that $V_0\cap\Omega$ is pseudoconvex. Let $\{\Omega_\nu\}_{\nu=1}^\infty$ be an exhaustion of $V_0\cap\Omega$ by strictly pseudoconvex domains with smooth boundary. Let $f = \sum_J' f_J d\bar{z}_J \in C^\infty_{(0,q)}(\overline{\Omega}_\nu)\cap\mathrm{dom}(\bar{\partial}^*_\nu)$, where we use a subscript $\nu$ to indicate operators on $\Omega_\nu$. A fundamental estimate due to Catlin ([5], Theorem 2.1, [4], inequality (2.3); see also [1], section 2 for a somewhat different approach to this type of estimate) gives for any function $\chi \in C^2(\overline{\Omega}_\nu)$, $0 \leq \chi \leq 1$:

$$\sum_K{}' \int_{\Omega_\nu} \sum_{j,k} \frac{\partial^2\chi}{\partial z_j \partial\bar{z}_k} f_{jK}\bar{f}_{kK} \leq C(\|\bar{\partial}f\|^2_{0,\Omega_\nu} + \|\bar{\partial}^*_\nu f\|^2_{0,\Omega_\nu}), \tag{12}$$

with a constant independent of $\nu$. The prime indicates summation over strictly increasing $(q-1)$ tuples, and we adopt the usual convention that the coefficients of $f$ be defined for all multi-indices so as to be antisymmetric. Regularizing the plurisubharmonic function $\lambda$ (from the hypothesis of the theorem) on $\Omega_\nu$ as in [12], proof of Theorem 4.4.2, and invoking (2) and (12) we find

$$\int_{\Omega_\nu} \delta^{-2\varepsilon}|f|^2 \leq C(\|\bar{\partial}f\|^2_{0,\Omega_\nu} + \|\bar{\partial}^*_\nu f\|^2_{0,\Omega_\nu}), \tag{13}$$

again with $C$ independent of $\nu$. If we could apply (13) to the restriction of $\varphi u$ to $\Omega_\nu$, the right-hand side of (11) would be seen to be dominated by $\|\bar{\partial}u\|_0 + \|\bar{\partial}^*u\|_0$, and we would be done. The problem is of course that while $\varphi u$ is smooth on $\overline{\Omega}_\nu$ (since the coefficients of $u$ are harmonic), it need not be in the domain of $\bar{\partial}^*_\nu$. To rectify the situation, we define forms $u_\nu$ on $\Omega_\nu$ as follows:

$$u_\nu := \bar{\partial}N_{q-1,\nu}\theta(\varphi u) + \bar{\partial}^*_\nu N_{q+1,\nu}\bar{\partial}(\varphi u). \tag{14}$$

When $q = 1$, we use the $\bar{\partial}$-Neumann operator $N_0$ on functions, see for example [16], Proposition 2.5. Then $u_\nu \in \mathrm{dom}(\bar{\partial}_\nu) \cap \mathrm{dom}(\bar{\partial}^*_\nu)$ and $u_\nu \in C^\infty_{(0,q)}(\overline{\Omega}_\nu)$ (since $\Omega_\nu$ is strictly pseudoconvex, the smoothness up to the boundary is preserved by $N_\nu$). Thus we can apply (13) to $u_\nu$. Taking into account that $\bar{\partial}_\nu^* \bar{\partial} N_{q-1,\nu}$ and $\bar{\partial}\bar{\partial}_\nu^* N_{q+1,\nu}$ are orthogonal projections in the respective $\mathcal{L}^2$-spaces (hence have norm 1), the result is

$$\int_{\Omega_\nu} \delta^{-2\varepsilon} |u_\nu|^2 \leq C(\|\bar{\partial}(\varphi u)\|_{0,\Omega_\nu}^2 + \|\theta(\varphi u)\|_{0,\Omega_\nu}^2) \tag{15}$$

$$\leq C(\|\bar{\partial}(\varphi u)\|_{0,V_0 \cap \Omega}^2 + \|\theta(\varphi u)\|_{0,V_0 \cap \Omega}^2).$$

We think of the forms $u_\nu$ as forms on $V_0 \cap \Omega$ by setting them zero on $(V_0 \cap \Omega) \setminus \Omega_\nu$. Then, because the $u_\nu$ are bounded in $\mathcal{L}^2_{(0,q)}(V_0 \cap \Omega)$ independently of $\nu$ (because $\bar{\partial} N_{q-1,\nu}$ and $\bar{\partial}_\nu^* N_{q+1,\nu}$ have bounds on their norms in the respective $\mathcal{L}^2$-spaces depending only on the diameter of $\Omega_\nu$), a suitable subsequence will converge weakly in $\mathcal{L}^2_{(0,q)}(V_0 \cap \Omega)$. Observing that $\varphi u \in \text{dom}(\bar{\partial}) \cap \text{dom}(\bar{\partial}^*)$ on $V_0 \cap \Omega$, one checks that this limit equals $\varphi u$ (as in [15], proof of Theorem 1; compare also [2], proof of part (3) of Theorem 2). Combining this with (15) gives

$$\|\delta^{-\varepsilon} \varphi u\|_{0,V_0 \cap \Omega}^2 \leq C(\|\bar{\partial} u\|_{0,V_0 \cap \Omega}^2 + \|\bar{\partial}^* u\|_{0,V_0 \cap \Omega}^2 + \|u\|_{0,V_0 \cap \Omega}^2). \tag{16}$$

In view of (11), this completes the proof of Theorem 1 (since again $\|u\|_{0,\Omega}$ is dominated by $\|\bar{\partial} u\|_{0,\Omega} + \|\bar{\partial}^* u\|_{0,\Omega}$).

## 4. Remaining proofs

To prove Theorem 2, it suffices to see that for $1 \leq k \leq m$, there exist bounded plurisubharmonic functions $\lambda_k$, defined near $P$ in $\Omega_k = \{\rho_k(z) < 0\}$ with $dd^c \lambda_k \geq \text{const.}(\delta_k(z))^{-2\varepsilon_k} dd^c |z|^2$ (as currents), where $\rho_k(z)$ denotes the distance to the boundary of $\Omega_k$. The function $\lambda := \lambda_1 + \cdots + \lambda_m$, defined near $P$ in $\Omega$, then has all the properties required in Theorem 1 (with $\varepsilon := \min_{1 \leq k \leq m} \varepsilon_k$). The existence of the functions $\lambda_k$ follows readily from Catlin's fundamental theorem ([5], Theorem 9.2) on good plurisubharmonic functions near finite type points in a smooth boundary. We fix $k$, but in order to avoid cluttering the notation, we temporarily suppress the subscript $k$. According to [5], Theorem 9.2, there exists $\varepsilon > 0$ and a neighborhood $V$ of $P$ such that for all sufficiently small $\delta > 0$ there is a smooth plurisubharmonic function $g_\delta$ in $V$ with $0 \leq |g_\delta| \leq 1$ and

$$dd^c g_\delta(z) \geq C \delta^{-2\varepsilon} dd^c |z|^2, \qquad z \in V \cap S_\delta. \tag{17}$$

Here $C$ is a constant independent of $\delta$, and $S_\delta$ denotes the set of points in $\Omega$ (near $P$) whose distance to the boundary is less than $\delta$. Choose $k_0 \in N$ such that $2^{-k_0}$ is sufficiently small. Let $0 < \varepsilon' < \varepsilon$. Set

$$g := \sum_{k=k_0}^{\infty} 2^{-2k(\varepsilon-\varepsilon')} g_{2^{-k}}. \tag{18}$$

Because $|g_{2^{-k}}| \leq 1$ on $V$, the sum in (18) converges uniformly on $V$, and $g$ is a bounded continuous plurisubharmonic function in $V$. Estimating $dd^c g$ from below by $dd^c(2^{-2k(\varepsilon-\varepsilon')} g_{2^{-k}})$ on the open sets $(S_{2^{-k}} \setminus \overline{S}_{2^{-k-2}}) \cap V$ (which cover $S_{2^{-k_0}} \cap V$) shows that, as currents on $S_{2^{-k_0}} \cap V$,

$$dd^c g \geq C(\delta(z))^{-2\varepsilon'} dd^c |z|^2. \tag{19}$$

Returning to the situation of Theorem 2, we have shown the existence of the plurisubharmonic functions $\lambda_k$, $1 \leq k \leq m$, with the required properties, and the proof of the theorem is complete.

Corollary 3 is a direct consequence of the proof of Theorem 1. Since $W_0^1(\Omega) \hookrightarrow \mathcal{L}^2(\Omega)$ compactly, we again only have to consider forms with harmonic coefficients. For a positive constant $M$, let $U$ be the neighborhood of the boundary given in the assumptions of Corollary 2. For a cut-off function $\varphi \in C_0^\infty(\Omega)$ identically one in a neighborhood of $\Omega \setminus U$, the arguments in the proof of Theorem 1 based on (12) and (14) give

$$M \int_\Omega |(1-\varphi)u|^2 \leq C(\|\bar{\partial}(1-\varphi)u\|_0^2 + \|\bar{\partial}^*(1-\varphi)u\|_0^2)$$
$$\leq C(\|\bar{\partial}u\|_0^2 + \|\bar{\partial}^*u\|_0^2) + C\|\nabla\varphi \cdot u\|_0^2. \tag{20}$$

Here $\nabla\varphi \cdot u$ denotes terms which are sums of (components of) $u$ times derivatives of $\varphi$. Note that $\operatorname{supp} \nabla\varphi \subseteq \operatorname{supp} \varphi$ which is compact in $\Omega$. Taking into account that the restriction to compact subsets is compact (in $\mathcal{L}^2$-norm) for harmonic functions, (20) shows that a sequence with $\|\bar{\partial}u_n\|_0 + \|\bar{\partial}^*u_n\|_0$ bounded has a subsequence that converges in $\mathcal{L}^2_{(0,q)}(\Omega)$. Consequently, the $\bar{\partial}$-Neumann operator $N_q$ is compact. It should be noted that although this argument is phrased slightly differently from Catlin's ([4], proof of Theorem 1), the only new ingredient is the "regularization" (14) which allows to derive (20) for nonsmooth $\Omega$ by passing to smooth subdomains.

## 5. Further remarks

*Remark 2:* Michel and Shaw ([15], Theorem 2) point out that their proof works in the case of strictly pseudoconvex domains with Lipschitz boundary, i.e. domains with a Lipschitz defining function whose Hessian is bounded from below by const. $dd^c|z|^2$.

*Remark 3:* The results in [11] are stated for the piecewise $C^2$ case, but the authors indicate how to get by with $C^1$. The condition on the singular part of the boundary is that its $(2n-1)$-dimensional Euclidean volume be zero. These results are for domains in a complex Hermitian manifold.

*Remark 4:* In the strictly pseudoconvex case, the relevant plurisubharmonic (defining) functions have Hessians bounded away from 0. By composing with $-(-x)^{1-2\varepsilon}$, one obtains plurisubharmonic functions as in Theorem 1 (and thus $\varepsilon$-subellipticity, when the boundary is Lipschitz) for all $\varepsilon < 1/2$. It should be possible to obtain $\varepsilon = 1/2$ essentially by observing that the measures $\delta^{-2\varepsilon}dV$, when *normalized*, approach surface measure on the boundary as $\varepsilon$ tends to zero. However, in this situation, it is more natural to exploit the boundary term in the Kohn-Morrey inequality in the first place, which is the approach taken in [15] and [11].

Department of Mathematics
Texas A&M University
College Station, TX 77843-3368
straube@math.tamu.edu